\documentclass[11pt, reqno]{amsart}
\usepackage{latexsym}
\usepackage{amssymb}
\usepackage{amsthm}
\usepackage{amsmath}
\usepackage{graphicx}
\usepackage{color}
\usepackage[all]{xy}
\allowdisplaybreaks   

\newtheorem{theorem}{Theorem}[section]
\newtheorem*{Theorem}{Theorem}
\newtheorem{corollary}[theorem]{Corollary}
\newtheorem*{Corollary}{Corollary}
\newtheorem{lemma}[theorem]{Lemma}

\newtheorem{example}[theorem]{Example}
\newtheorem{definition}[theorem]{Definition}

\newtheorem{proposition}[theorem]{Proposition}

\newtheorem*{Remark}{Remark}

\DeclareMathOperator{\Hom}{Hom}
\DeclareMathOperator{\Tor}{Tor}
\DeclareMathOperator{\Ext}{Ext}
\DeclareMathOperator{\Ker}{Ker}
\DeclareMathOperator{\Imm}{Im}
\DeclareMathOperator{\HH}{HH}
\DeclareMathOperator{\Ho}{H}
\newcommand{\bk}{\mathbf{k}}
\newcommand{\kG}{\bk G}
\newcommand{\D}{\mathcal{D}}
\newcommand{\ba}{\mathbf{a}}

\begin{document}
\begin{abstract}

\small{Let $A$ be any finite dimensional Hopf algebra over a field $\bk$. We specify the Tate and Tate-Hochschild cohomology for $A$ and introduce cup products that make them become graded rings. We establish the relationship between these rings. In particular, the Tate-Hochschild cohomology of $A$ is isomorphic (as algebras) to its Tate cohomology with coefficients in an adjoint module. Consequently, the Tate cohomology ring of $A$ is a direct summand of its Tate-Hochschild cohomology ring. As an example, we explicitly compute both the Tate and Tate-Hochschild cohomology for the Sweedler algebra $H_4$.  
}

\end{abstract}

\title[Tate and Tate-Hochschild Cohomology]
{Tate and Tate-Hochschild Cohomology\\
for finite dimensional Hopf Algebras}

\author{Van C. Nguyen}
\address{Department of Mathematics\\Texas A\&M University \\College Station, TX 77843-3368}
\email{vcnguyen@math.tamu.edu}

\maketitle

\section{Introduction}
\label{sec:intro}

Tate cohomology was introduced by John Tate in 1952 for group cohomology arising from class field theory \cite{Tate}. Others then generalized his theory to the group ring $RG$ where $R$ is a commutative ring and $G$ is a finite group. Unlike the usual cohomology, this theory is based on complete resolutions, and hence, yields cohomology groups in both positive and negative degrees. Over the past several decades, a great deal of effort has gone into the study of this new cohomology. A summary may be found in (\cite{Brown}, Ch.~VI) or (\cite{CaEi}, Ch.~XII). In the early 1980's, through an unpublished work, Pierre Vogel extended the Tate cohomology to any group and even to any ring using unbounded chain complexes. For finite groups, the Tate-Vogel cohomology coincides with the Tate cohomology. Accounts of Vogel's construction appeared, for examples, in a paper by Goichot \cite{Goi} and in another paper by Benson and Carlson \cite{BeCa} in 1992. In the 1980's, Buchweitz introduced another construction of Tate cohomology of a two-sided Noetherian and Gorenstein ring, using the stable module category influenced by the work of Auslander and Bridger (\cite{Buch}, Section~6). Many authors have also considered the Hochschild analogue of Tate cohomology for Frobenius algebras. For instance, one of the first attempts was given in Nakayama's paper in 1957 on the complete cohomology of Frobenius algebras, using the complete standard complex (or complete bar resolution) \cite{Nakayama}. The stable Hochschild cohomology of a Frobenius algebra, using the stable module category, was studied in various papers, e.g. \cite{EuSc}. More recently, using complete resolutions, Bergh and Jorgensen defined the Tate-Hochschild cohomology of an algebra $A$ whose enveloping algebra $A^e$ is two-sided Noetherian and Gorenstein over a field $\bk$ \cite{BeJo}. If the Gorenstein dimension of $A^e$ is $0$, then this cohomology agrees with the usual Hochschild cohomology in positive degrees. It is noted in \cite{BeJo} that this Tate-Hochschild definition is equivalent to that using the stable module category in \cite{EuSc}, at least in the finite dimensional case. In this paper, we specialize the Tate and Tate-Hochschild cohomology to finite dimensional Hopf algebras over a field $\bk$ based on the definition from \cite{BeJo}. Moreover, since any finite dimensional Hopf algebra is a Frobenius algebra (\cite{Mont}, Theorem~2.1.3), results from \cite{BeJo}, \cite{EuSc}, and \cite{Nakayama} apply.

We begin Section~\ref{sec:prelim} by introducing some definitions and notation. In Section~\ref{sec:Tate}, we generalize the construction in (\cite{Ben2}, Section~5.15) to define the Tate cohomology $\widehat{\Ho}^n(A,\bk)$ of a finite dimensional Hopf algebra $A$:
 \[\widehat{\Ho}^n(A,\bk):= \widehat{\Ext}_A^n(\bk,\bk), \text{ for all } n \in \mathbb{Z}, \]
where the $\Ext$ functor is taken over an $A$-complete resolution of $\bk$. We observe that the Tate cohomology groups share some common properties as those of the usual cohomology groups. 

Let $A^e=A \otimes_\bk A^{op}$ be the enveloping algebra of $A$. We consider $A$ as an $A^e$-module and obtain its $A^e$-complete resolution. Based on the definition in \cite{BeJo}, we define the $n$-th Tate-Hochschild cohomology group of $A$ as:
 \[\widehat{\HH}^n(A,A) := \widehat{\Ext}^n_{A^e}(A,A), \text{ for all } n \in \mathbb{Z}. \]
 
It is known that the usual cohomology of $A$ is a direct summand of its Hochschild cohomology (\cite{GiKu}, Prop.~5.6 and Cor.~5.6). This motivates us to ask if the same assertion holds for the Tate cohomology version. Section~\ref{sec:Tate-Hochschild} is where we start making the first comparison between the Tate and the Tate-Hochschild cohomology groups of $A$. In particular, by Proposition~4.2, $\widehat{\Ho}^n(A,\bk)$ is a vector space direct summand of $\widehat{\HH}^n(A,A)$, for all $n>0$. On the other hand, by (\cite{BeJo}, Cor.~3.8), if $\nu$ is the Nakayama automorphism of a Frobenius algebra $A$ such that $\nu^2=\textbf{1}$, then there is an isomorphism:
 \[\widehat{\HH}^n(A,A) \cong \widehat{\HH}^{-(n+1)}(A,A), \text{ for all } n \in \mathbb{Z}. \]
Following such a symmetric result, we further study the Nakayama automorphism in Section~\ref{sec:Nakayama}. We compute ``formulas'' for the Nakayama automorphisms of the Sweedler algebra $H_4$ and of the general Taft algebra, in an attempt to relate the Tate-Hochschild cohomology groups in positive degrees to those in negative degrees. This computation does not lead us to conclude the cohomology relation for a general finite dimensional Hopf algebra. However, we find that it is helpful in verifying our results for the Tate-Hochschild cohomology of the Sweedler algebra $H_4$ in Section~\ref{sec:Taft algebra}, as in this case we have $\nu^2=\textbf{1}$.

We approach a broader setting by establishing cup products for the two Tate cohomology types in Section~\ref{sec:products}. These multiplication structures turn $\widehat{\Ho}^*(A,\bk)$ and $\widehat{\HH}^*(A,A)$ into graded rings. Using the ring structures, we prove the following isomorphism in Section~\ref{sec:relationship}:
\begin{Theorem}
Let $A$ be a finite dimensional Hopf algebra over a field $\bk$. Let $A^{ad}$ be $A$ as an $A$-module via the left adjoint action. Then there exists an isomorphism of algebras:
 \[ \widehat{\HH}^*(A,A) \cong \widehat{\Ho}^*(A,A^{ad}). \]
\end{Theorem}
\vspace{-0.5em}
The same assertion holds in the usual cohomology for finite dimensional Hopf algebras (\cite{GiKu}, Prop.~5.6), in particular, it is well-known in finite group cohomology (\cite{EiMac}, I.5.7). This theorem generates an important relation between the Tate and the Tate-Hochschild cohomology rings:
\begin{Corollary}
If $A$ is a finite dimensional Hopf algebra over a field $\bk$, then $\widehat{\Ho}^*(A,\bk)$ is a direct summand of $\widehat{\HH}^*(A,A)$ as a module over $\widehat{\Ho}^*(A,\bk)$.
\end{Corollary}
\vspace{-0.5em}
Therefore, the Tate and Tate-Hochschild cohomology share the same relation as that of the usual cohomology. In the last section, we compute the Tate and Tate-Hochschild cohomology for the Taft algebra, in particular, the Sweedler algebra $H_4$. These examples help us to understand the above relation better, and once again, to illustrate the symmetric result in \cite{BeJo}.

\section{Preliminaries and Notations}
\label{sec:prelim}

Projective resolutions are commonly used to compute the cohomology of an algebra. Here, we apply a more general resolution, which involves both positive and negative degrees:
\begin{definition}
\label{complete resolution}
Suppose $R$ is a two-sided Noetherian ring. A \textbf{complete resolution} of a finitely generated $R$-module $M$ is an exact complex $\mathbb{P}=\{\{P_i\}_{i \in \mathbb{Z}}, d_i: P_i \rightarrow P_{i-1}\}$ of finitely generated projective $R$-modules such that:
\begin{enumerate}
 \item The dual complex $\Hom_R(\mathbb{P},R)$ is also exact.
 \item There exists a projective resolution $\mathbb{Q} \stackrel{\varepsilon}{\rightarrow} M$ of $M$ and a chain map $\mathbb{P} \stackrel{\varphi}{\rightarrow} \mathbb{Q}$ where $\varphi_n$ is bijective for $n \geq 0$ and $0$ for $n<0$.
\end{enumerate}
\end{definition}

\vspace{-0.5em}
Unlike projective resolutions, complete resolutions, in general, do not always exist. However, given the setting as in Definition~\ref{complete resolution}, Theorems 3.1 and 3.2 in \cite{AvMa} guarantee the existence of such complete resolutions. 

In \cite{Tate}, Tate introduced a cohomology theory that can be defined by using complete resolutions as follows (\cite{Ben2}, Section~5.15). Let $G$ be a finite group and $R$ be a commutative ring with $G$ acts trivially on $R$. If
 \[\cdots \xrightarrow{d_3} P_2 \xrightarrow{d_2} P_1 \xrightarrow{d_1} P_0 \xrightarrow{\varepsilon} R \rightarrow 0 \]
is an $RG$-projective resolution of $R$, then apply $\Hom_R(-,R)$ to get a dual sequence:
 \[0 \rightarrow R \rightarrow \Hom_R(P_0,R) \rightarrow \Hom_R(P_1,R) \rightarrow \Hom_R(P_2,R) \rightarrow \cdots \] 
This is again an exact sequence of projective $RG$-modules. Splicing these two sequences together, one forms a doubly infinite sequence:
 \[\cdots \xrightarrow{d_3} P_2 \xrightarrow{d_2} P_1 \xrightarrow{d_1} P_0 \rightarrow \Hom_R(P_0,R) \rightarrow \Hom_R(P_1,R) \rightarrow \cdots\]
Introducing the notation $P_{-(n+1)} := \Hom_R(P_n,R)$, we arrive at a complete resolution of $R$:
 \[\mathbb{P}: \hspace{1 cm} \cdots \xrightarrow{d_3} P_2 \xrightarrow{d_2} P_1 \xrightarrow{d_1} P_0 \rightarrow P_{-1} \rightarrow P_{-2} \rightarrow \cdots\] 
 
Fixing a (left) $RG$-module $M$ and applying $\Hom_{RG}(-,M)$ to $\mathbb{P}$ produces a new complex. The homology of this new complex is the \textbf{Tate cohomology} for $RG$: 
 \[\widehat{\Ho}^n(G,M) := \widehat{\Ext}^n_{RG}(R,M) = \Ho^n(\Hom_{RG}(\mathbb{P},M)), \text{ for all } n \in \mathbb{Z}.\]

As previously discussed, the goal of this paper is to study this cohomology theory for finite dimensional Hopf algebras. For the rest of this paper, we let $\bk$ be our base field. Tensor products will be over $\bk$ unless stated otherwise. Let $A$ be a finite dimensional Hopf algebra over $\bk$ with counit $\varepsilon: A \rightarrow \bk$, and coproduct $\Delta: A \rightarrow A \otimes A$, expressed in Sweedler notation as $\displaystyle \Delta(a) = \sum a_1 \otimes a_2$. Here, $\bk$ is a trivial (left) $A$-module via the counit map, that is, if $a\in A, r \in \bk$, then $a \cdot r :=\varepsilon(a)r$. The antipode of $A$, denoted by $S$, is bijective in the finite dimensional case (\cite{Mont}, Theorem~2.1.3), and its inverse is denoted by $\overline{S}$. The $\bk$-dual $\Hom_\bk(-,\bk)$ is denoted by $D(-)$ and the ring-dual $\Hom_A(-,A)$ is denoted by $(-)^*$. This is an unfortunate notation, because to a Hopf algebraist, $D(A)$ usually stands for the ``Drinfeld double'' of $A$. However, to be consistent with the primary reference paper \cite{BeJo}, we will adopt the above notation. We observe that $D(A)$ is also a finite dimensional Hopf algebra, where the algebra structure of $A$ becomes the coalgebra structure of $D(A)$, the antipode of $A$ translates into an antipode $D(S)$ of $D(A)$ in a canonical fashion, and so on. Moreover, for any left $A$-module $M$, the dual module $D(M)$ is a left $A$-module with the action of $A$ given by the antipode, that is, $(a \cdot f)(m)= f(S(a)m)$, where $a \in A$, $m \in M$, and $f \in D(M)$. $A$ is a Frobenius algebra (\cite{Mont}, Theorem~2.1.3), which is self-injective, so projective $A$-modules and injective $A$-modules coincide.

The opposite algebra $A^{op}$ has the same underlying set and linear operation as $A$ but with multiplication performed in the reverse order: $a*b = ba$, for all $a$ and $b \in A$. Let $A^e = A \otimes A^{op}$ denote the enveloping algebra of $A$ and define $\sigma: A \rightarrow A^e$ by $\sigma(a) = \sum a_1 \otimes S(a_2)$. Checking that $\sigma$ is an injective algebra homomorphism, we may identify $A$ with the subalgebra $\sigma(A)$ of $A^e$. Moreover, we can induce $A^e$-modules from $A$-modules as follows. Let $M$ be a left $A$-module and consider $A^e$ as a right $A$-module via right multiplication by $\sigma(A)$. Then $A^e \otimes_A M$ is a left $A^e$-module, with $A^e$-action given by $a \cdot(b \otimes_A m) = ab \otimes_A m, \text{ for all } a, b \in A^e$, and $m \in M$.

We will use the following notation for the usual (non-Tate) cohomology and Hochschild cohomology of the Hopf algebra $A$, respectively: 
 \begin{align*}
 \Ho^*(A,M)&:= \Ext^*_A(\bk,M)=\displaystyle \bigoplus_{n\geq0} \Ext^n_A(\bk,M), \\
 \HH^*(A,M)&:= \Ext^*_{A^e}(A,M)=\displaystyle \bigoplus_{n\geq0} \Ext^n_{A^e}(A,M),
 \end{align*}
where $M$ denotes a left $A$-module in the former case and an $A$-bimodule in the latter case.

\section{Tate Cohomology for finite dimensional Hopf Algebras}
\label{sec:Tate}

\subsection{Definition}
\label{subsec:Tate construction}

Generalizing the construction in (\cite{Ben2}, Section~5.15), we can explicitly form a complete resolution $\mathbb{P}$ of $\bk$ as an $A$-module from an $A$-projective resolution $P_{\bullet}$ of $\bk$. This can be done by splicing $P_{\bullet}$ with its dual complex $D(P_{\bullet})$, which is also an exact sequence of finitely generated projective $A$-modules, since $A$ is a self-injective algebra. One can check that the resulting complex $\mathbb{P}$ is exact and satisfies the definition of a complete resolution of $\bk$. We define the Tate cohomology groups of $A$ with coefficients in a left $A$-module $M$ as:
 \[\widehat{\Ho}^n(A,M):= \widehat{\Ext}_A^n(\bk,M) =\Ho^n(\Hom_A(\mathbb{P},M)), \text{ for all } n \in \mathbb{Z}. \]
The Tate homology groups $\widehat{\Ho}_n(A,M):=\widehat{\Tor}_n^A(\bk,M)$ are defined analogously by applying $- \otimes_A M$ to $\mathbb{P}$ and taking the $n$-th homology of the new complex. Here, we are only interested in the Tate cohomology. Observe that in our context, naturally, the Tate (co)homology does not depend on the choice of the projective resolution of $\bk$ (by the ordinary Comparison Theorem), and hence, is independent of the complete resolution of $\bk$ (\cite{AvMa}, Theorem~5.2 and Lemma~5.3). 

\begin{Remark}{\em
Instead of using complete resolutions, there is another formulation of the Tate cohomology for $A$ via the stable module category (\cite{Buch} Lemma~6.1.2). If $M$ and $N$ are finitely generated $A$-modules, we define $\underline{\Hom}_A(N,M)$ to be the quotient of $\Hom_A(N,M)$ by homomorphisms that factor through a projective module. Then for any integer $n$:
 \[ \widehat{\Ext}_A^n(\bk,M) \cong \underline{\Hom}_A(\Omega^n\bk,M) \cong \underline{\Hom}_A(\bk,\Omega^{-n}M), \]
or equivalently (\cite{Buch}, Prop.~6.5.1)
 \[ \widehat{\Ext}_A^n(\bk,M) \cong \lim_{\stackrel{\longrightarrow}{k,\; k+n \geq 0}} \underline{\Hom}_A(\Omega^{k+n}\bk,\Omega^k M), \] 
where $\Omega$ is the Heller operator, mapping an $A$-module to the kernel of a projective cover of that module. This definition, which is equivalent to that using complete resolutions, is useful especially when proving some results on the cochain level. However, to be consistent with \cite{BeJo} and to develop an alternative approach that may prove to be useful, we will use the definition with complete resolutions.}
\end{Remark}
\vspace{-1em}

\subsection{Properties of Tate Cohomology}
\label{subsec:prop}

We compare the Tate cohomology and the usual cohomology of $A$. Here are some important remarks:
 \begin{enumerate}
 \renewcommand{\labelenumi}{(\alph{enumi})}  
  \item For all $n>0$, $\widehat{\Ho}^n(A,M) \cong \Ho^n(A,M)$.
  
  \item The group $\widehat{\Ho}^0(A,M)$ is a quotient of $\Ho^0(A,M)$.
    
   These follow from the construction of complete resolutions.
    
  \item For all $n < -1$, we have isomorphisms: $\widehat{\Ho}^n(A,M) \cong H_{-(n+1)}(A,M)$ (\cite{Brown}, Prop.~I.8.3c). 
  
     \item If $0 \rightarrow M \rightarrow M' \rightarrow M'' \rightarrow 0$ is a short exact sequence of (left) $A$-modules, then there is a doubly infinite long exact sequence of Tate cohomology groups (\cite{AvMa}, Prop.~5.4) or (see \cite{Nakayama}, Theorem~1 for the Tate-Hochschild version):
 \[\cdots \rightarrow \widehat{\Ho}^n(A,M) \rightarrow \widehat{\Ho}^n(A,M') \rightarrow \widehat{\Ho}^n(A,M'') \rightarrow \widehat{\Ho}^{n+1}(A,M) \rightarrow \cdots \] 

 \item If $(N_j)_{j\in J}$ is a finite family of (left) $A$-modules and $(M_i)_{i\in I}$ is any family of $A$-modules, then there are natural isomorphisms, for all $n \in \mathbb{Z}$ (\cite{AvMa}, Prop.~5.7):
  \[\displaystyle \widehat{\Ext}^n_A(\bigoplus_{j \in J}{N_j}, M) \cong \prod_{j \in J} \widehat{\Ext}^n_A(N_j,M), \]
  \[\displaystyle \widehat{\Ext}^n_A(N, \prod_{i \in I} M_i) \cong \prod_{i \in I} \widehat{\Ext}^n_A(N,M_i). \] 
 \end{enumerate}

\section{Tate-Hochschild Cohomology for finite dimensional Hopf Algebras}
\label{sec:Tate-Hochschild}

\begin{definition}
Let $\bk, A$, and $A^e$ be defined as before (such that $A^e$ is two-sided Noetherian and Gorenstein). Let $M$ be an $A$-bimodule. For any integer $n \in \mathbb{Z}$, the $n$-th Tate-Hochschild cohomology group is defined as:
 \[\widehat{\HH}^n(A,M) := \widehat{\Ext}^n_{A^e}(A,M), \text{ for all } n \in \mathbb{Z}. \]
\end{definition}

\vspace{-0.5em}
When $A^e$ is a two-sided Noetherian and Gorenstein ring, by (\cite{AvMa}, Theorems 3.1, 3.2), every finitely generated $A^e$-module admits a complete resolution. Hence, we obtain a complete resolution $\mathbb{X}$ for $A$ as an $A^e$-module. Moreover, any bimodule $M$ of $A$ can be viewed as a left $A^e$-module by setting $(a \otimes b) \cdot m = amb$, for $a \otimes b \in A^e$ and $m \in M$. The $n$-th Tate-Hochschild cohomology group is the $n$-th homology group of the complex $\Hom_{A^e}(\mathbb{X},M)$. As noted in \cite{BeJo}, the assumption that $A^e$ is two-sided Noetherian is not necessary in the finite dimensional case. This also intrigues us to focus on finite dimensional Hopf algebras.

We recall the fact that a finite dimensional Hopf algebra $A$ is Frobenius which is self-injective. Lemmas~3.1 and~3.2 in \cite{BeJo} show that $A^{op}, A^{e}$ are also Frobenius, hence self-injective.  By the definition of Gorenstein ring, we see that finite dimensional self-injective algebras $A, A^{op}, A^{e}$ are Gorenstein of Gorenstein dimension 0. Therefore, the Tate-Hochschild cohomology groups of $A$ agree with the usual Hochschild cohomology groups in all positive degrees \cite{BeJo}:
 \[\widehat{\HH}^n(A,M) \cong \HH^n(A,M), \text{ for all } n>0. \]
\begin{Remark}{\em
The Tate analog of the Hochschild cohomology of a Frobenius algebra is also considered in (\cite{Nakayama}, Section~3) using complete standard complex, or in (\cite{EuSc}, 2.1.11) as the stable Hochschild cohomology $\underline{\Hom}_{A^e}(\Omega^n A,M)$. Eu and Schedler also showed that cup product and contraction structures extend to the stable, $\mathbb{Z}$-graded setting for the Tate-Hochschild cohomology ring of $A$ (\cite{EuSc}, Theorem~2.1.15).}
\end{Remark} 

\vspace{-0.5em}
We make our first attempt to relate the Tate and Tate-Hochschild cohomology:
\begin{proposition}
Let $A$ be a finite dimensional Hopf algebra over a field $\bk$. Then $\widehat{\Ho}^n(A,\bk)$ is a vector space direct summand of $\widehat{\HH}^n(A,A), \text{ for all } n>0$. That means, we may view $\widehat{\Ho}^n(A,\bk)$ as (isomorphic to) a quotient of $\widehat{\HH}^n(A,A)$, for all $n>0$.
\end{proposition}

\vspace{-1em}
\begin{proof}
For $n>0, \widehat{\Ho}^n(A,\bk) \cong \Ho^n(A,\bk)$ by \ref{subsec:prop} (a), and $\widehat{\HH}^n(A,A) \cong \HH^n(A,A)$ by the above discussion. Moreover, for a finite dimensional Hopf algebra $A$, its antipode $S$ is bijective. Applying (\cite{PW}, Lemma~7.2), we identify $\Ho^*(A,\bk)$ with a subalgebra of $\HH^*(A,A)$, and use appropriate isomorphism maps to obtain the desired result.
\end{proof}

\vspace{-0.5em}
We can simplify this relation even more. By (\cite{BeJo}, Cor.~3.8), if $\nu$ is the Nakayama automorphism of a Frobenius algebra $A$ such that $\nu^2=\textbf{1}$, then there is an isomorphism:
 \[\widehat{\HH}^n(A,A) \cong \widehat{\HH}^{-(n+1)}(A,A), \text{ for all } n \in \mathbb{Z}. \]
Thus, giving this symmetric structure, for $n>0, \widehat{\Ho}^n(A,\bk)$ is a vector space direct summand of $\widehat{\HH}^n(A,A) \cong \widehat{\HH}^{-(n+1)}(A,A)$. This leads us to examine when the Nakayama automorphism $\nu$ of a finite dimensional Hopf algebra has the property $\nu^2=\textbf{1}$. This is trivial if $A$ is symmetric, as $\nu$ is the identity in this case. We then assume $A$ is not symmetric.

\section{Nakayama Automorphism}
\label{sec:Nakayama}

Recall that for any Hopf algebra $A$, its dual $D(A)$ obtains the structure of $A$-module via a left action $a \rightharpoondown f$ and a right action $f \leftharpoondown a$ of $A$ on $D(A)$ as follows. For $a,b \in A, f \in D(A)$, the actions are given by: $(a \rightharpoondown f)(b) := f(S(a)b)$ and $(f \leftharpoondown a)(b) := f(bS(a)) = (S(a) \rightharpoonup f)(b)$. Therefore, $D(A)$ is a left and right $A$-module. There are also a natural left action $f \rightharpoonup a := \sum f(a_2)a_1$ and a right action $a \leftharpoonup f := \sum f(a_1)a_2$ of $D(A)$ on $A$. 

The counit $\varepsilon: A \rightarrow \bk$ gives rise to defining the space of \textbf{left integrals} of $A$: 
\[\int^l_A := \{x \in A \;|\;hx=\varepsilon(h)x, \text{ for all } h \in A \},\]
and the space of \textbf{right integrals}: 
\[\int^r_A := \{x \in A\;|\;xh=\varepsilon(h)x, \text{ for all } h \in A \}.\]
When $A$ is a finite dimensional Hopf algebra, we have $\text{dim}_\bk \int^l_A = \text{dim}_\bk \int^r_A = \text{dim}_\bk \int^l_{D(A)} = \text{dim}_\bk \int^r_{D(A)} = 1$, (\cite{Mont}, Theorem~2.1.3). This shows the existence of non-zero integrals in any finite dimensional Hopf algebra $A$ and also in its dual $D(A)$. 

Let $0 \neq f \in \int^l_{D(A)}$ be a non-zero left integral of $D(A)$, the non-degenerate associative bilinear Frobenius form $\mathcal{B}: A \times A \rightarrow \bk$ of $A$ is given by: $\mathcal{B}(x,b) = f(xb), \text{ for all } x, b \in A$. The \textbf{Nakayama automorphism} $\nu: A \rightarrow A$ satisfies $\mathcal{B}(x,b) = \mathcal{B}(b,\nu(x))$. Replacing $\mathcal{B}$ with a new Frobenius form $\mathcal{B}'$ defined by a unit element $u \in A$ gives us a new automorphism $\nu' = I_u \circ \nu$, where $I_u$ is the inner automorphism $r \mapsto uru^{-1}$. The Nakayama automorphism $\nu$ is unique up to composition with an inner automorphism. Equivalently, it is a well-defined element of the group of outer automorphisms of $A$. 

Let $t \in \int^r_A$, we have $at \in \int^r_A, \text{ for all } a \in A$. The (right) \textbf{modular function} for $A$ is an algebra homomorphism $\alpha \in D(A)$ such that $at=\alpha(a)t, \text{ for all } a \in A$. By (\cite{FiMoSc}, Lemma~1.5), the Nakayama automorphism $\nu: A \rightarrow A$ of a finite dimensional Hopf algebra has the form: 
 \[\nu(a) = \overline{S}^2(a \leftharpoonup \alpha) = (\overline{S}^2 a) \leftharpoonup \alpha, \]
and its inverse is defined by 
 \[\nu^{-1}(a) = S^2(a \leftharpoonup \alpha^{-1}) = (S^2a) \leftharpoonup \alpha^{-1}. \]
This formula for $\nu$ also shows that $\nu$ has a finite order dividing $2\cdot\text{dim}_\bk(A)$.

Consider for every $a\in A$:
 \[\nu^2(a) = \overline{S}^2((\overline{S}^2 a) \leftharpoonup \alpha) \leftharpoonup \alpha = (\overline{S}^4 a) \leftharpoonup \alpha^2. \]
Therefore, to get $\nu^2=\textbf{1}$, we need the antipode $S$ to be of order 4 (since $S$ is bijective, $\overline{S}$ is also of order 4), and the right modular function $\alpha$ to be of order 2. Keeping in mind that these are sufficient but not necessary conditions for $\nu^2=\textbf{1}$, we examine several known cases for the former condition:
\begin{itemize}
 \item If $A$ is commutative or cocommutative then $S^2=\textbf{1}$, (\cite{Mont}, Cor.~1.5.12).
 \item If $A$ is a finite-dimensional semisimple and cosemisimple Hopf algebra over a field $\bk$ of characteristic 0 or of characteristic $p > (\text{\text{dim}}_\bk A)^2$, then $S^2=\textbf{1}$, (\cite{Mont}, Theorem~2.5.3).
\end{itemize}
 
\begin{example}[Sweedler algebra]
\label{ex:Sweedler algebra}
{\em  
The smallest non-commutative, non-cocommutative Hopf algebra with dimension 4 over a given field $\bk$ of characteristic $\neq 2$ was described by Sweedler as:
 \[H_4 = \bk \left\langle 1,g,x,gx\;|\;g^2=1, x^2=0, xg=-gx \right\rangle \]
with the Hopf algebra structure:
\begin{align*}
\Delta(g) &= g\otimes g, 
&\Delta(x)&=1\otimes x + x\otimes g, \\
\varepsilon(g)&=1, 
&\varepsilon(x)&=0, \\
S(g)&=g=g^{-1}, 
&S(x)&=-xg.
\end{align*}
The antipode $S$ has order 4. The integral spaces of $H_4$ are:
 \[\int^l_{H_4} = \bk\langle x+gx\rangle, \hspace{2cm} \int^r_{H_4} = \bk\langle x-gx\rangle. \]
To find the right modular function for $H_4$, we take the right integral $t=x-gx$ and multiply this element by all the basis elements of $H_4$ on the left to see what the right modular function looks like. The desired modular function $\alpha \in D(H_4)$ should satisfy: $at=\alpha(a)t, \text{ for all } a \in H_4$. After checking for each basis element, we have:
 \[\alpha(x) = \alpha(gx) = \alpha(xg) = 0, \hspace{2cm} \alpha(g)=-1, \hspace{2cm} \alpha(1)=1,\]
and \vspace{-0.5em}
 \[\nu(g) = -g, \hspace{2cm} \nu(x) = -x. \]
By direct computation or by using the fact that $S$ has order 4, we obtain:
\[\nu^2(a)= (\overline{S}^4 a) \leftharpoonup \alpha^2 = a, \hspace{1cm} \text{for all } a \in H_4. \]
Therefore, $\nu^2=\textbf{1}$ on $H_4$. It follows from the above discussion that the Tate-Hochschild cohomology of $H_4$ is symmetric in the sense that:
 \[\widehat{\HH}^n(H_4,H_4) \cong \widehat{\HH}^{-(n+1)}(H_4,H_4), \text{ for all } n \in \mathbb{Z}. \]
An explicit computation of the Tate-Hochschild cohomology of $H_4$ in Section~\ref{sec:Taft algebra} will illustrate this result again.}
\end{example}

\vspace{-0.5em}
We examine a more general case, the Taft algebra, as described in (\cite{Mont2}, Example~2.9):

\begin{example}
\label{ex:Taft algebra}
{\em
Let $N \geq 2$ be a positive integer. Assume the field $\bk$ contains a primitive $N$-th root of unity $\omega$. Let $A$ be the Taft algebra generated over $\bk$ by two elements $g$ and $x$, subject to the following relations:
 \[A = \bk \left\langle 1, g, x, gx\;|\;g^N=1, x^N=0, xg=\omega gx \right\rangle. \]
$A$ is a Hopf algebra with structure given by: 
 \begin{align*}
  \Delta(g)&= g \otimes g,
  &\Delta(x)&= 1 \otimes x + x \otimes g, \\
  \varepsilon(g)&= 1,
  &\varepsilon(x)&= 0,\\
  S(g)&= g^{-1},
  &S(x)&=-xg^{-1}.
 \end{align*}
Note that $A$ is not semisimple and is of dimension $N^2$. Using $\displaystyle t=\sum^{N-1}_{j=0} x^{N-1}g^j$ as a right integral of $A$, we can check that the right modular function $\alpha \in D(A)$ is defined as: 
 \[\alpha(x) = \alpha(gx) = \alpha(xg) = 0, \hspace{2cm} \alpha(g) = \omega, \hspace{2cm} \alpha(1)=1. \]
The Nakayama map $\nu: A \rightarrow A$ is given by: 
 \[\nu(g) = \omega g, \hspace{2cm} \nu(x) = \omega x, \]
which yields: \vspace{-0.5em}
$$\nu^2(a)= \begin{cases} 
            1 &\quad a=1\\
            \omega^2a &\quad a=g,x\\
            \omega^4a &\quad a=xg.
            \end{cases}$$ 
So $\nu^2$ is not the identity map unless $N=2$, which we have seen in the previous example. We remark that under some other settings (for example, using $\overline{S}$ instead of $S$, or using the left integral and left modular function), one might obtain different ``formulas'' for the Nakayama automorphism, such as $\nu'(g) = \omega g$ and $\nu'(x) = \omega^{-1}x$ as described in (\cite{Mont2}, Example~2.9).}
\end{example}

\vspace{-0.5em}
It is not easy to classify all non-symmetric finite dimensional Hopf algebras $A$ such that $\nu^2=\textbf{1}$. What we have done here does not effectively yield a general connection between the Tate and Tate-Hochschild cohomology of a finite dimensional Hopf algebra. In the next section, we observe that these two Tate-cohomology types obtain ring structures which can help us to develop a deeper understanding of their relation as algebras.

\section{Cup Products}
\label{sec:products}
\subsection{Cup Product on Tate Cohomology}
\label{subsec:Tate product}

Suppose $\mathbb{P}$ is an $A$-complete resolution of $\bk$. Based on the discussion in (\cite{Brown}, Section~VI.5), we also note the following difficulties in constructing the cup product on Tate cohomology:

First of all, $\mathbb{P} \otimes \mathbb{P}$ is not a complete resolution of $\bk \otimes \bk \cong \bk$, as $(\mathbb{P} \otimes \mathbb{P})_{+}$ is not the same as the tensor product of resolutions $\mathbb{P}_{+} \otimes \mathbb{P}_{+}$, where $\mathbb{P}_{+} = \{P_n\}_{n\geq0}$. Consequently, using the map $\Hom_A(\mathbb{P},M) \otimes \Hom_A(\mathbb{P},N) \rightarrow \Hom_A(\mathbb{P} \otimes \mathbb{P}, M \otimes N)$ would not obviously induce a cohomology product in Tate cohomology as it does in the usual non-Tate cohomology. Secondly, when applying the diagonal approximation (a chain map that preserves augmentation) $\Gamma: \mathbb{P} \rightarrow \mathbb{P} \otimes \mathbb{P}$, for any $n \in \mathbb{Z}$, there are infinitely many $(i,j)$ such that $i+j=n$, and the dimension-shifting property in Section~\ref{subsec:prop} suggests that the corresponding cup products should all be non-trivial. So $\Gamma$ should have a non-trivial component $\Gamma_{ij}, \text{ for all } (i,j)$. Hence, the range of $\Gamma$ should be the graded module which is $\prod_{i+j=n} P_i \otimes P_j$ in the dimension $n$, rather than $\bigoplus_{i+j=n} P_i \otimes P_j$. This discussion motivates us to the following definitions:

Let $\varepsilon: \mathbb{P} \rightarrow \bk$ be an $A$-complete resolution of $\bk$ and let $d$ be the differentials in $\mathbb{P}$. We form the \textbf{complete tensor product} $\mathbb{P} \widehat{\otimes} \mathbb{P}$ by defining:
 \[(\mathbb{P} \widehat{\otimes} \mathbb{P})_n = \displaystyle \prod_{i+j=n} P_i \otimes P_j, \text{ for all } n \in \mathbb{Z}, \]
with the ``total differential'' $\partial_{i,j} = d^v_{i,j} + d^h_{i,j}$, where $d^h_{i,j}=d_i \widehat{\otimes} \textbf{1}_\mathbb{P}$ and $d^v_{i,j}=(-1)^i\textbf{1}_\mathbb{P} \widehat{\otimes} d_j$. It can be easily seen that $\mathbb{P} \widehat{\otimes} \mathbb{P}$ is an acyclic complex of $A$-modules.

On the other hand, given graded modules $B, B', C, C'$ and module homomorphisms $u: C \rightarrow B$ of degree $r$ and $v: C' \rightarrow B'$ of degree $s$, there is a map $u \widehat{\otimes} v: C \widehat{\otimes} C' \rightarrow B \widehat{\otimes} B'$ of degree $r+s$ defined by: \vspace{-0.5em}
 \[ (u \widehat{\otimes} v)_n = \displaystyle \prod_{i+j=n} (-1)^{is} u_i \otimes v_j: \displaystyle \prod_{i+j=n} C_i \otimes C'_j \rightarrow \displaystyle \prod_{i+j=n} B_{i+r} \otimes B'_{j+s}. \]

\begin{definition}
A \textbf{complete diagonal approximation map} is a chain map $\Gamma: \mathbb{P} \rightarrow \mathbb{P} \widehat{\otimes} \mathbb{P}$ such that $(\varepsilon \widehat{\otimes} \varepsilon) \circ \Gamma_{0} = \varepsilon$, that is, $\Gamma$ is an augmentation-preserving chain map.
\end{definition}

\vspace{-0.5em}
A similar argument to the proof given in (\cite{Brown}, Section~VI.5) shows the existence of such complete diagonal approximation map $\Gamma: \mathbb{P} \rightarrow \mathbb{P} \widehat{\otimes} \mathbb{P}$. Let $M$ and $N$ be left $A$-modules. Then $M \otimes N$ is also a left $A$-module via the coproduct: $a \cdot (m \otimes n)= \sum{a_1m \otimes a_2n}$, for all $a \in A, m \in M$, and $n \in N$. We define a cochain cup product:
 \[\smile: \Hom_A(P_i,M) \otimes \Hom_A(P_j,N) \rightarrow \Hom_A(P_{i+j},M \otimes N) \] 
given by 
 \[f \smile g = (f \widehat{\otimes} g) \circ \Gamma, \]
with $f \in \Hom_A(P_i,M)$ and $g \in \Hom_A(P_j,N)$ represent elements of $\widehat{\Ho}^i(A,M)$ and $\widehat{\Ho}^j(A,N)$, respectively. One verifies that by the definition of differentials on the total complex, the usual coboundary formula holds:
 \[\delta(f \smile g) = (\delta f) \smile g + (-1)^i f \smile (\delta g). \]
It follows from the formula that the product of two cocycles is again a cocycle and the product of a cocycle with a coboundary is a coboundary. Thus, this induces a well-defined product on Tate cohomology $\widehat{\Ho}^i(A,M) \otimes \widehat{\Ho}^j(A,N) \rightarrow \widehat{\Ho}^{i+j}(A,M \otimes N)$. Moreover, this cup product is unique, in the sense that it is compatible with connecting homomorphisms in the long exact cohomology sequence, independent of the choice of $\mathbb{P}$ and $\Gamma$, associative, that is, 
\[(f \smile g) \smile h = f \smile (g \smile h),\]
and $1 \in \widehat{\Ho}^0(A,\bk)$ is an identity. One proves that using the dimension-shifting property in Section~\ref{subsec:prop} and exactness of tensor products over $\bk$, similarly as in (\cite{Brown},~V.3.3 and Lemma~VI.5.8). It is immediate from the definitions that this product is natural with respect to coefficient homomorphisms. For example, an $A$-module homomorphism $M \otimes N \rightarrow Q$ yields products
\[\widehat{\Ho}^i(A,M) \otimes \widehat{\Ho}^j(A,N) \rightarrow \widehat{\Ho}^{i+j}(A,Q)\]
by compositing the cup product and the induced map $\widehat{\Ho}^{i+j}(A,M \otimes N) \rightarrow \widehat{\Ho}^{i+j}(A,Q)$. In particular, when $M=N=\bk$, $\widehat{\Ho}^*(A,\bk)$ is a graded ring; and when $N=\bk$, $\widehat{\Ho}^*(A,M)$ is a graded module over $\widehat{\Ho}^*(A,\bk)$.

\subsection{Cup Product on Tate-Hochschild Cohomology}
\label{subsec:Tate-Hochschild product}

There is also a cup product on the Tate-Hochschild cohomology:
\[\widehat{\HH}^i(A,M) \otimes \widehat{\HH}^j(A,N) \rightarrow \widehat{\HH}^{i+j}(A,M \otimes_A N)\]
where $M$ and $N$ are $A$-bimodules (which can be viewed as (left) $A^e$-modules). Before showing this cup product, let us recall some useful lemmas whose proofs can be found in \cite{PW}. We provide the sketch of the proofs for the sake of completeness. Let $\sigma: A \rightarrow A^e$ be defined by $\sigma(a) = \sum a_1 \otimes S(a_2)$. Recall from Section~\ref{sec:prelim} that $A^e$ may be viewed as a right $A$-module via right multiplication by elements of $\sigma(A)$. 
\begin{lemma}
\label{Tate lemma 1}
$A \cong A^e \otimes_A \bk$ as left $A^e$-modules, where $A^e \otimes_A \bk$ is the induced $A^e$-module.
\end{lemma}

\vspace{-1em}
\begin{proof}
Let $f: A \rightarrow A^e \otimes_A \bk$ be the function defined by:
 \[f(a)=a \otimes 1 \otimes_A 1,\]
and let $g: A^e \otimes_A \bk \rightarrow A$ be the function defined by:
 \[g(a \otimes b \otimes_A 1) = ab,\]
for all $a, b \in A$. One can easily check that $f$ and $g$ are both $A^e$-module homomorphisms, and that they are inverses of each other.
\end{proof}

\begin{lemma}
\label{Tate lemma 2}
$A^e$ is a (right) projective $A$-module.
\end{lemma}

\vspace{-1em}
\begin{proof}
Since $A$ is finite dimensional, its antipode map $S$ is bijective. Moreover, $S$ is an $A$-module map: for all $a, b \in A, S(ab) = S(b)S(a)=S(a)*S(b)$ in $A^{op}$. This implies $S: A \rightarrow A^{op}$ is an isomorphism of right $A$-modules, where $A$ acts on $A$ by right by multiplication and on $A^{op}$ by multiplication by $S(A)$. This yields an isomorphism of right $A$-modules: $A \otimes A \rightarrow A \otimes A^{op} = A^e$. Since $A$ is projective over itself, and the tensor product with a projective module is projective (\cite{Ben1}, Prop.~3.1.5), $A \otimes A$ is a projective right $A$-module. Therefore, $A^e$ is a projective right $A$-module, where $A$ acts on $A^e$ by multiplying $\sigma(A)$.
\end{proof}

\vspace{-0.5em}
Let $\mathbb{X}$ be any $A^e$-complete resolution of $A$. By the same argument as in Section~\ref{subsec:Tate product}, $\mathbb{X} \widehat{\otimes}_A \mathbb{X}$ is an acyclic chain complex of $A^e$-modules. Since $\otimes_A$ is not an exact functor in general, the existence of a complete diagonal approximation map $\Gamma$ does not follow trivially from \cite{Brown} as before. We show it here in detail.
\begin{lemma}
\label{Tate lemma 3}
There exists a complete diagonal approximation map $\Gamma: \mathbb{X} \rightarrow \mathbb{X} \widehat{\otimes}_A \mathbb{X}$.
\end{lemma}

\vspace{-1em}
\begin{proof}
Observe that $A^e \otimes_A A^e = (A \otimes A^{op}) \otimes_A (A \otimes A^{op}) = A \otimes (A^{op} \otimes_A A) \otimes A^{op} \cong A \otimes A^{op} \otimes A^{op} \cong A^e \otimes_\bk A$. Since $A^e$ acts only on the outermost two factors of $A$, dropping $A$ in the third step does not change the $A^e$-module structure. Therefore, $A^e \otimes_A A^e \cong  A^e \otimes_\bk A$ is an $A^e$-module isomorphism, not just a $\bk$-module isomorphism. 

As $A$ is a free (hence projective) $\bk$-module, $A^e \otimes_\bk A$ is also free as an $A^e$-module. Consequently, $A^e \otimes_A A^e$ is a free $A^e$-module. In general, tensor product over $A$ of free $A^e$-modules is free. Since any projective module is a direct summand of a free module, this implies that for all $i, j \in \mathbb{Z}$, $X_i \otimes_A X_j$ is projective as an $A^e$-module. Again, because $A^e$ is self-injective, projective $A^e$-modules are also injective. Therefore, $X_i \otimes_A X_j$ is injective, implying the direct product $(\mathbb{X} \widehat{\otimes}_A \mathbb{X})_n$ is an injective $A^e$-module for all $n \in \mathbb{Z}$.

As remarked in (\cite{Ben1}, Theorem~2.4.2), to form the chain map $\Gamma_+: \mathbb{X}_+ \rightarrow (\mathbb{X} \widehat{\otimes}_A \mathbb{X})_+$ in non-negative degrees, it suffices for the complex $\mathbb{X}_+$ to consist of projective modules but it need not be exact, and for the complex $(\mathbb{X} \widehat{\otimes}_A \mathbb{X})_+$ to be exact but not necessarily to consist of projective modules. Since $\mathbb{X}_+$ is a projective resolution of $A$, we can apply the ordinary Comparison Theorem to obtain a chain map $\Gamma_+$ that is augmentation-preserving. We then consider the projective $A^e$-modules in negative degrees of these complexes, which are (relatively) injective as discussed above. By a generalization of (\cite{Brown}, Prop.~VI.2.4), the family of maps $\Gamma_+$ extends to a complete chain map $\Gamma: \mathbb{X} \rightarrow \mathbb{X} \widehat{\otimes}_A \mathbb{X}$ in both positive and negative degrees. 
\end{proof}

\vspace{-0.5em}
We may define a cup product on Tate-Hochschild cohomology as follows. Let $M$ and $N$ be $A$-bimodules, then $M \otimes_A N$ is also an $A$-bimodule which can be considered as a left $A^e$-module via $(a \otimes b) \cdot (m \otimes_A n) = am \otimes_A nb$, for $a \otimes b \in A^e, m \in M$ and $n \in N$. Let $f \in \Hom_{A^e}(X_i,M)$ represent an element of $\widehat{\HH}^i(A,M)$ and let $g \in \Hom_{A^e}(X_j,N)$ represent an element of $\widehat{\HH}^j(A,N)$. Then:
 \[ f \smile g = (f \widehat{\otimes} g) \circ \Gamma \in \Hom_{A^e}(X_{i+j},M \otimes_A N) \] represents an element of $\widehat{\HH}^{i+j}(A,M \otimes_A N)$. One can check that this product is independent of $\mathbb{X}$ and $\Gamma$ and satisfies certain properties as in Section~\ref{subsec:Tate product}. In particular, if $M=N=A$, then $\widehat{\HH}^*(A,A)$ is a graded ring, for $A \cong A \otimes_A A$. If $N=A$, then $\widehat{\HH}^*(A,M)$ is a graded $\widehat{\HH}^*(A,A)$-module.

\section{Relationship between Tate and Tate-Hochschild Cohomology of $A$}
\label{sec:relationship}

We begin with a lemma based on the original Eckmann-Shapiro Lemma but generalized to a complete resolution:

\begin{lemma}[Eckmann-Shapiro]
\label{Tate lemma 4}
Let $B$ be a ring, let $C \subseteq B$ be a subring for which $B$ is flat as a right $C$-module. Let $M$ be a left $C$-module and let $N$ be a left $B$-module. Consider $N$ to be a left $C$-module via restriction of the action, and let $B \otimes_C M$ denote the induced $B$-module where $B$ acts on the leftmost factor by multiplication. Then for all $i \in \mathbb{Z}$, there is an isomorphism of abelian groups:
 \[\widehat{\Ext}^i_C(M,N) \cong \widehat{\Ext}^i_B(B \otimes_C M,N). \]
If $B$ and $C$ are $\bk$-algebras, then this is an isomorphism of vector spaces over $\bk$. 
\end{lemma}

\vspace{-1em}
\begin{proof}
Let $\varepsilon: \mathbb{P} \rightarrow M$ be a $C$-complete resolution of $M$. Since $B \otimes_C C \cong B$ as a left $B$-module, the induced modules $B \otimes_C P_i$ are projective $B$-modules, for all $i \in \mathbb{Z}$. The induced complex $B \otimes_C \mathbb{P}$ is exact as $B$ is flat over $C$, with the ``augmentation map'' $\textbf{1}_B \otimes_C \varepsilon: B \otimes_C \mathbb{P} \rightarrow B \otimes_C M$. So it is a complete resolution of $B \otimes_C M$ as a $B$-module.

It suffices to show that for all $i \in \mathbb{Z}$, $\Hom_C(P_i,N) \cong \Hom_B(B \otimes_C P_i,N)$ as abelian groups. This follows from the Nakayama relations (\cite{Ben1}, Prop.~2.8.3). One should also check that these isomorphisms commute with the differentials. By the definition of Tate cohomology, these isomorphisms will comprise a chain map that induces an isomorphism on cohomology and give us the desired result.
\end{proof}

\vspace{-0.5em}
We consider $A$ to be an $A$-module by the left adjoint action: for $a, b \in A$, $a \cdot b = \sum a_1bS(a_2)$, and denote this $A$-module by $A^{ad}$. More generally, if $M$ is an $A$-bimodule, denote by $M^{ad}$ the left $A$-module with action given by $a \cdot m = \sum a_1mS(a_2)$, for all $a \in A, m \in M$. We now prove our main result:

\begin{theorem}
\label{Tate iso}
Let $A$ be a finite dimensional Hopf algebra over a field $\bk$. Then there exists an isomorphism of algebras:
 \[ \widehat{\HH}^*(A,A) \cong \widehat{\Ho}^*(A,A^{ad}). \]
\end{theorem}

\vspace{-1em}
\begin{proof}
By Lemma~\ref{Tate lemma 2}, $A^e$ is a projective, hence, flat $A$-module. We then apply Lemmas \ref{Tate lemma 1} and \ref{Tate lemma 4} with $B=A^e$, $C=A$ is identified as a subalgebra of $A^e$, $M=\bk$ is a left $A$-module, and $N=A \cong A^e \otimes_A \bk$ is an induced left $A^e$-module. We get $\widehat{\Ext}^*_A(\bk,A^{ad}) \cong \widehat{\Ext}^*_{A^e}(A^e \otimes_A \bk,A)$ as $\bk$-modules, i.e. $\widehat{\Ho}^*(A,A^{ad}) \cong \widehat{\HH}^*(A,A)$ as $\bk$-modules.

To show this is an algebra isomorphism, it remains to prove that cup products are preserved by this isomorphism. Let $\mathbb{P}$ denote an $A$-complete resolution of $\bk$. Since $A^e$ is a (right) projective $A$-module by Lemma~\ref{Tate lemma 2}, $\mathbb{X} = A^e \otimes_A \mathbb{P}$ is an $A^e$-complete resolution of $A^e \otimes_A \bk \cong A$.

Recall that $A$ is acting on $A^e$ on the left as well as on the right via $\sigma$. We define an $A$-chain map $\iota: \mathbb{P} \rightarrow \mathbb{X}$ by $\iota(p)= (1 \otimes 1) \otimes_A p, \text{ for all } p \in P_i, i \in \mathbb{Z}$. Let $f \in \Hom_{A^e}(X_i,A)$ be a cocycle representing a cohomology class in $\widehat{\Ext}^i_{A^e}(A,A)$. The corresponding cohomology class in $\widehat{\Ext}^i_{A}(\bk,A^{ad})$ is represented by $f \circ \iota$.

Let $\Gamma: \mathbb{P} \rightarrow \mathbb{P} \widehat{\otimes} \mathbb{P}$ be a complete diagonal approximation map. $\Gamma$ induces a cup product on cohomology as discussed in Section~\ref{subsec:Tate product}. $\Gamma$ also induces a chain map $\Gamma': \mathbb{X} \rightarrow \mathbb{X} \widehat{\otimes}_A \mathbb{X}$ as follows. There is a map of $A^e$-chain complexes $\phi: A^e \otimes_A (\mathbb{P} \widehat{\otimes} \mathbb{P}) \rightarrow \mathbb{X} \widehat{\otimes}_A \mathbb{X}$ given by: 
 \[\phi((a \otimes b) \otimes_A (p \otimes q)) = ((a \otimes 1) \otimes_A p) \otimes_A ((1 \otimes b) \otimes_A q). \]
$\Gamma$ induces a map from $A^e \otimes_A \mathbb{P}$ to $A^e \otimes_A (\mathbb{P} \widehat{\otimes} \mathbb{P})$. Let $\Gamma'$ be the composition of this map with $\phi$. 

Let $f \in \Hom_{A^e}(X_i,A)$ and $g \in \Hom_{A^e}(X_j,A)$ be cocycles. The above discussions imply the following diagram commutes:
\begin{displaymath} 
  \xymatrix{\mathbb{X} \ar[r]^{\Gamma' \quad}               & \mathbb{X} \widehat{\otimes}_A \mathbb{X} \ar[r]^{f \widehat{\otimes} g}           & A \otimes_A A \ar[r]^{\quad \sim} & A \\
            \mathbb{P} \ar[r]^{\Gamma \quad} \ar[u]^{\iota} & \mathbb{P} \widehat{\otimes} \mathbb{P} \ar[r]^{(f\iota) \widehat{\otimes} (g\iota)} & A \otimes A \ar[r]^{\quad m}        & A \ar@{=}[u]}
\end{displaymath}
where $m$ denotes the multiplication $a \otimes b \stackrel{m}{\longmapsto} ab, \text{ for all } a, b \in A$. 

As described in Section~\ref{sec:products}, the top row yields a product in $\widehat{\Ext}^*_{A^e}(A,A)$ and the bottom row yields a product in $\widehat{\Ext}^*_A(\bk,A^{ad})$. Thus, cup products are preserved and $\widehat{\HH}^*(A,A)$ is isomorphic to $\widehat{\Ho}^*(A,A^{ad})$ as algebras.
\end{proof}

\vspace{-0.5em}
As a consequence of Theorem~\ref{Tate iso}, to determine the Tate-Hochschild cohomology of a finite dimensional Hopf algebra $A$, it suffices to compute its Tate cohomology with coefficients in the adjoint $A$-module. One may, therefore, apply known examples for Tate cohomology groups (such as, \cite{CaEi}, Section~XII.7) to compute the corresponding Tate-Hochschild cohomology. Furthermore, we arrive at the desired relation between the Tate and Tate-Hochschild cohomology rings of $A$: 

\begin{corollary}
If $A$ is a finite dimensional Hopf algebra over a field $\bk$, then $\widehat{\Ho}^*(A,\bk)$ is a direct summand of $\widehat{\HH}^*(A,A)$ as a module over $\widehat{\Ho}^*(A,\bk)$.
\end{corollary}

\vspace{-1em}
\begin{proof}
Under the left adjoint action of $A$ on itself, the trivial module $\bk$ is isomorphic to the submodule of $A^{ad}$ given by all scalar multiples of the identity 1. In fact, $\bk$ is a direct summand of $A^{ad}$ with complement the augmentation ideal $\Ker(\varepsilon)$, where $\varepsilon: A \rightarrow \bk$ is the counit map. 

$\widehat{\Ext}^*_A(\bk,-)$ is additive. By Theorem~\ref{Tate iso}, we have:
\begin{align*}
\widehat{\HH}^*(A,A) \cong \widehat{\Ho}^*(A,A^{ad}) &= \widehat{\Ext}^*_A(\bk,A^{ad}) \\
                                                  &\cong \widehat{\Ext}^*_A(\bk,\bk) \oplus \widehat{\Ext}^*_A(\bk,\Ker(\varepsilon)) \\
                                                  &\cong \widehat{\Ho}^*(A,\bk) \oplus \widehat{\Ext}^*_A(\bk,\Ker(\varepsilon)).
\end{align*}
Both $\bk$ and $\Ker(\varepsilon)$ are closed under multiplication, and this multiplication induces multiplications on $\widehat{\Ext}^*_A(\bk,-)$ which is also compatible with the ring structure on $\widehat{\Ext}^*_A(\bk,A^{ad})$. Hence, this is in fact a direct summand, where $\widehat{\HH}^*(A,A) \cong \widehat{\Ho}^*(A,A^{ad})$ is considered as a (left) module over $\widehat{\Ho}^*(A,\bk)$ with action via (left) multiplication.                                 
\end{proof}

\vspace{-0.5em}
For the rest of this paper, we will compute the Tate and Tate-Hochschild cohomology for the Taft algebra, in particular, the Sweedler algebra. These simple examples should give the readers a rough procedure to produce other examples.

\section{Example: Taft Algebra}
\label{sec:Taft algebra}
\subsection{Tate Cohomology}
\label{subsec:Taft Tate}

We provide an example to compute the Tate and Tate-Hochschild cohomology groups of a Taft algebra $A$, which is a Hopf algebra of dimension $N^2$:
\[A = \bk \left\langle 1, g, x, gx\;|\;g^N=1, x^N=0, xg=\omega gx \right\rangle \]
as described in Example~\ref{ex:Taft algebra}, where the field $\bk$ contains a primitive $N$-th root of unity $\omega$. 
It is known that as an algebra, a Taft algebra is a smash product $A = B\#\kG$ (or a skew group algebra), with $B=\bk[x]/(x^N)$, and $G$ is a finite cyclic group generated by $g$ of order $N$ acting on $B$. Let $\chi: G \rightarrow \bk^\times$ be the character, that is, a group homomorphism from $G$ to the multiplicative group of $\bk$, defined by $\chi(g) = \omega$. $G$ acts by automorphisms on $B$ via:
\[ ^gx= \chi(g)x  = \omega x.\]
Note that since $G$ is generated by $g$, all $G$-actions reduce to the action of the generator $g$. By the definition of smash product, $A$ is $B \otimes \kG$ as a vector space, with the multiplication:
\[ (b_1\otimes g_1)(b_2 \otimes g_2) = b_1(^{g_1}b_2) \otimes g_1g_2, \]
for all $b_1,b_2 \in B$ and $g_1, g_2 \in G$. We abbreviate $b_i \otimes g_i$ by $b_ig_i$. Moreover, as the characteristic of $\bk$ does not divide $|G|=N$,
$\kG$ is semisimple and so all the cohomology of $\kG$ is trivial except in the degree 0. From (\cite{Ste}, Cor.~3.4), the Taft algebra's usual cohomology is known as:
 \[\displaystyle \Ho^*(A,\bk) \cong (\Ext^*_B(\bk,\bk))^{G},\]
where the superscript $G$ denotes the invariants under the action of $G$. Again, as the characteristic of $\bk$ is relatively prime to $|G|$, $G$-invariants may be taken in a complex prior to taking the cohomology. We consider the following $B$-free resolution of $\bk$:
\begin{equation}
\label{eq:Tate5}
 \cdots \xrightarrow{\cdot x} B \xrightarrow{\cdot x^{N-1}} B \xrightarrow{\cdot x} B \xrightarrow{\cdot x^{N-1}} B \xrightarrow{\cdot x} B \xrightarrow{\varepsilon} \bk \rightarrow 0,
\end{equation}
with $\varepsilon(x)=0$ is the augmentation map. This resolution could be extended to a projective $A$-resolution of $\bk$ by giving $B$ the following actions of $G$, for all $b \in B$ and $i>0$: 
\begin{itemize}
 \item In degree $0$, $g \cdot b := (^gb)$.
 \item In degree $2i$, $g \cdot b :=\chi(g)^{iN}(^gb) = (^gb)$, since $\chi(g)^{iN} = \omega^{iN} = 1^i=1$.
 \item In degree $2i+1$, $g \cdot b :=\chi(g)^{iN+1}(^gb) = \omega(^gb)$.
\end{itemize}
One checks that this group action commutes with the differentials in (\ref{eq:Tate5}) in each degree. Thus, we may extend the differentials $\cdot x^{N-1}$ and $\cdot x$ in (\ref{eq:Tate5}) to be $A$-module maps. Moreover, since the characteristic of $\bk$ does not divide the order of $G$, an $A$-module is projective if and only if its restriction to $B$ is projective. With these actions, (\ref{eq:Tate5}) is indeed an $A$-projective resolution of $\bk$.

Take the $\bk$-dual $\Hom_\bk(-,\bk)$ of (\ref{eq:Tate5}), we have:
\begin{equation}
\label{eq:Tate6}
 0 \rightarrow \bk \xrightarrow{D(\varepsilon)} D(B) \xrightarrow{D(\cdot x)} D(B) \xrightarrow{D(\cdot x^{N-1})} D(B) \xrightarrow{D(\cdot x)} D(B) \xrightarrow{D(\cdot x^{N-1})} \cdots
\end{equation}
which is an exact sequence of projective $A$-modules. Splicing (\ref{eq:Tate5}) and (\ref{eq:Tate6}) together, we obtain an $A$-complete resolution of $\bk$:
\begin{equation}
\label{eq:Tate7}
  \mathbb{P}: \hspace{0.5cm} \cdots \xrightarrow{\cdot x} B \xrightarrow{\cdot x^{N-1}} B \xrightarrow{\cdot x} B \xrightarrow{\xi} D(B) \xrightarrow{D(\cdot x)} D(B) \xrightarrow{D(\cdot x^{N-1})} D(B) \xrightarrow{D(\cdot x)} \cdots 
\end{equation}
where $\xi = D(\varepsilon) \circ \varepsilon$. To compute the Tate cohomology of $A$ with coefficients in $\bk$, apply $\Hom_A(-,\bk)=\overline{(-)}$ to (\ref{eq:Tate7}):
\begin{equation}
\label{eq:Tate8}
 \cdots \xrightarrow{\overline{D(\cdot x)}} \overline{D(B)} \xrightarrow{\overline{D(\cdot x^{N-1})}} \overline{D(B)} \xrightarrow{\overline{D(\cdot x)}} \overline{D(B)} \xrightarrow{\overline{\xi}} \overline{B} \xrightarrow{\overline{(\cdot x)}} \overline{B} \xrightarrow{\overline{(\cdot x^{N-1})}} \overline{B} \xrightarrow{\overline{(\cdot x)}} \cdots
\end{equation}
which is a complex of $A$-modules. 
Take the homology of this new complex, we will obtain the Tate cohomology of the Taft algebra $A$.

Let us compute explicitly for the case $N=2$. Here, $G \cong \mathbb{Z}_2$, $B=\bk[x]/(x^2)$, and $A = B\# \kG$ is the Sweedler algebra $H_4$ that we have seen in Example~\ref{ex:Sweedler algebra}. $B$ has a basis $\{1,x\}$ and $D(B)$ has a dual basis $\{f_1, f_x\}$. By previously defined actions of $G$, $B$ is an $H_4$-module with:
\begin{align*}
g\cdot x &=-x,   \hspace{.5cm} g\cdot 1 =1, \hspace{.5cm} \text{in 0 and even degrees} \\
g\cdot x &=x,    \hspace{.5cm} g\cdot 1 =-1, \hspace{.5cm} \text{in odd degrees}
\end{align*}
such that this action commutes with the differential maps in (\ref{eq:Tate5}). 
We denote $(-)_{ev}$ for objects in even degrees and $(-)_{odd}$ for objects in odd degrees, given the corresponding actions of $G$ as $H_4$-modules; hence, (\ref{eq:Tate7}) becomes:
\[\mathbb{P}: \hspace{0.5cm} \cdots \xrightarrow{\cdot x} B_{ev} \xrightarrow{\cdot x} B_{odd} \xrightarrow{\cdot x} B_0 \xrightarrow{\xi} D(B)_{odd} \xrightarrow{D(\cdot x)} D(B)_{ev} \xrightarrow{D(\cdot x)} D(B)_{odd} \xrightarrow{D(\cdot x)} \cdots  \]
as an $H_4$-complete resolution of $\bk$. $D(B)$ is an $H_4$-module via the action: $(g \cdot f)(b)=f(S(g)\cdot b)=f(g \cdot b)$, for $f \in D(B), b \in B, g \in G$ and $S(g) = g^{-1} = g$ in $H_4$. Checking on the basis elements, we see that the $G$-actions on $D(B)$ can be carried along: 
\begin{align*}
g\cdot f_1 &=f_1, 	\hspace{.5cm} g\cdot f_x=-f_x, \hspace{.5cm} \text{in 0 and even degrees} \\
g\cdot f_1 &=-f_1, 	\hspace{.5cm} g\cdot f_x=f_x, \hspace{.5cm} \text{in odd degrees}.
\end{align*}
By identifying $f_1 \leftrightarrow x$ and $f_x \leftrightarrow 1$, we have $D(B)_{ev} \cong B_{odd}$ and $D(B)_{odd} \cong B_{ev}$ as $H_4$-modules. As a result, $\mathbb{P}$ can be written as:
\[\mathbb{P}: \hspace{0.5cm} \cdots \xrightarrow{\cdot x} B_{ev} \xrightarrow{\cdot x} B_{odd} \xrightarrow{\cdot x} B_0 \xrightarrow{\xi} B_{ev} \xrightarrow{D(\cdot x)} B_{odd} \xrightarrow{D(\cdot x)} B_{ev} \xrightarrow{D(\cdot x)} \cdots  \]
Apply $\Hom_{H_4}(-,\bk)=\overline{(-)}$ to $\mathbb{P}$, we get the complex of $H_4$-modules:
\[\cdots \xrightarrow{\overline{D(\cdot x)}} \overline{B_{ev}} \xrightarrow{\overline{D(\cdot x)}} \overline{B_{odd}} \xrightarrow{\overline{D(\cdot x)}} \overline{B_{ev}} \xrightarrow{\overline{\xi}} \overline{B_0} \xrightarrow{\overline{(\cdot x)}} \overline{B_{odd}} \xrightarrow{\overline{(\cdot x)}} \overline{B_{ev}} \xrightarrow{\overline{(\cdot x)}} \cdots\]
For all $f \in \Hom_B(B,\bk)$ and $b \in B$, $1=f(b)=b\cdot f(1)=\varepsilon(b)f(1)$. We may identify $f$ with a map in $\Hom_\bk(\bk,\bk)$ and obtain an isomorphism $\Hom_B(B,\bk) \cong \Hom_\bk(\bk,\bk)$ which is isomorphic to $\bk$. We then observe that $\overline{B}=\Hom_{H_4}(B,\bk)$ is contained in $\Hom_B(B,\bk) \cong \bk$. One can check that under the corresponding group actions:
$$\overline{B} = \begin{cases}
                 \bk &\quad \text{ in 0 and even degrees} \\
                 0 &\quad \text{ in odd degrees.}\\
                 \end{cases} $$
This simplifies the above complex to:
\[\cdots \xrightarrow{\overline{D(\cdot x)}} \bk \xrightarrow{\overline{D(\cdot x)}} 0 \xrightarrow{\overline{D(\cdot x)}} \bk \xrightarrow{\overline{\xi}} \bk \xrightarrow{\overline{(\cdot x)}} 0 \xrightarrow{\overline{(\cdot x)}} \bk \xrightarrow{\overline{(\cdot x)}} \cdots\]
To compute the homology of this complex, we need to see what $\overline{\xi}$ looks like. $\overline{\xi}: \overline{D(B)} \rightarrow \overline{B}$ is defined as $\overline{\xi}(h)(b)=h(\xi(b))$ for $h\in \overline{D(B)}, b \in B$. By exactness of $\mathbb{P}$, $\Imm(\xi)=\Ker(D(\cdot x))=\left\langle f_1\right\rangle$. As an $H_4$-module map, $h$ sends $f_1 \mapsto 0$, and $f_x \mapsto 1$. Therefore, $\overline{\xi}(h)(b)=h(\xi(b))=h(f_1)=0$ and $\overline{\xi}$ is a 0-map.  

Put these together, we have computed the Tate cohomology for Sweedler algebra $H_4$:
$$\widehat{\Ho}^n(H_4,\bk) = \begin{cases}
                         \bk &\quad n<-1, n \text{ is odd}\\
                         0 &\quad n<-1, n \text{ is even}\\
                         \bk &\quad n =-1, 0\\
                         0 &\quad n>0, n \text{ is odd}\\
                         \bk &\quad n>0, n \text{ is even.}
                         \end{cases} $$
In other words, $\widehat{\Ho}^n(H_4,\bk) \cong \widehat{\Ho}^{-(n+1)}(H_4,\bk)$, for all $n \in \mathbb{Z}$.

\subsection{Tate-Hochschild Cohomology}
\label{subsec:Taft Tate-Hochschild}

In order to compute the Tate-Hochschild cohomology of a general Taft algebra $A$, we use the following subalgebra of $A^e = A \otimes A^{op}$ as in \cite{BW}: 
\[\D := B^e\#\kG \cong \bigoplus_{g\in G} (Bg \otimes Bg^{-1}) \subset A^e, \]
where the action of $G$ on $B^e$ is diagonal, that is, $^g(a\otimes b) = (^ga) \otimes (^gb)$. This isomorphism is given by $(a \otimes b)g \mapsto ag \otimes (^{g^{-1}}b)g^{-1}$, for all $a,b \in B$, and $g \in G$. Note that $B$ is a $\D$-module under left and right multiplications. Since the characteristic of $\bk$ does not divide $|G|=N$, the Hochschild cohomology $\HH^*(A) := \Ext^*_{A^e}(A,A)$ is known to satisfy:
\[\HH^*(A) \cong \Ext^*_{\D}(B,A) \cong \Ext^*_{B^e}(B,A)^G \]
as graded algebras. $\Ext^*_{B^e}(B,A)^G$ consists of invariants under the actions induced from the actions of $G$ on $\D$-modules, see (\cite{BW},~4.9) or (\cite{Ste}, Cor.~3.4) for more details. Observe the following $B^e$-free resolution of $B$, (\cite{Weibel}, Exercise~9.1.4):
\begin{equation}
\label{eq:Tate9}
 \cdots \rightarrow B^e \xrightarrow{\cdot v} B^e \xrightarrow{\cdot u} B^e \xrightarrow{\cdot v} B^e \xrightarrow{\cdot u} B^e \xrightarrow{m} B \rightarrow 0, 
\end{equation}
where $m$ is the multiplication map $a \otimes b \mapsto ab$, $$u = x \otimes 1 - 1 \otimes x, \qquad \qquad \qquad \text{and } v = x^{N-1} \otimes 1 + x^{N-2} \otimes x + \cdots + 1 \otimes x^{N-1}.$$ 

Using this resolution, one computes $\HH^n(B) \cong B/(x^{N-1})$ and $(\HH^n(B))^G \cong \bk$ for all $n>0$. This resolution also becomes a $\D$-projective resolution of $B$ by giving the following actions of $G$ on $B^e$, for all $a \otimes b \in B^e, g \in G$ and integers $i>0$:
\begin{itemize}
 \item In degree $0$, $g \cdot (a \otimes b) = (^ga) \otimes (^gb)$.
 \item In degree $2i$, $g \cdot (a \otimes b) = \chi(g)^{iN}(^ga) \otimes (^gb) = (^ga) \otimes (^gb)$, since $\chi(g)^{iN} = \omega^{iN} = 1^i= 1$.
 \item In degree $2i+1$, $g \cdot (a \otimes b) = \chi(g)^{iN+1}(^ga) \otimes (^gb) = \omega(^ga) \otimes (^gb)$.
\end{itemize} 
The differentials $\cdot u$ and $\cdot v$ commute with the group actions and turn out to be the maps of $\D$-modules. Since char($\bk$) does not divide $|G|$, a $\D$-module is projective if and only if its restriction to $B^e$ is projective. With these actions, (\ref{eq:Tate9}) becomes a $\D$-projective resolution of $B$. 

Because $B^e$ is a left $B$-module by multiplying by the leftmost factor in $B^e$, we can take the dual $\Hom_B(-,B)=(-)^*$ of (\ref{eq:Tate9}). Since $B \cong\Hom_B(B,B)$, we have:
\begin{equation}
\label{eq:Tate10}
 0 \rightarrow B \xrightarrow{m^*} (B^e)^* \xrightarrow{(\cdot u)^*} (B^e)^* \xrightarrow{(\cdot v)^*} (B^e)^* \xrightarrow{(\cdot u)^*} (B^e)^* \xrightarrow{(\cdot v)^*} \cdots
\end{equation}
One can show that this is an exact sequence of projective $\D$-modules. $\Hom_B(B^e,B)=(B^e)^*$ is a left $B^e$-module via action: 
\[ ((a \otimes b) \cdot f)(c \otimes d)= f((a \otimes b)(c \otimes d)) = f(ac \otimes db),\]
for all $f \in \Hom_B(B^e,B)$, and $a, b, c, d \in B$. The differentials $(\cdot u)^*, (\cdot v)^*,m^*$ are $B^e$-module homomorphisms, since they are just the compositions of maps, $d^*(f)=f \circ d$. 

For all $f \in \Hom_B(B^e,B)$, $f (a \otimes b) = f(a(1\otimes b)) = af(1 \otimes b)$. Hence, we may identify $1 \otimes b$ with $b \in B^{op}$ and $B^e$ with $\Hom_\bk(B^{op},B) \cong \Hom_B(B^e,B) = (B^e)^*$ under the correspondence $a \otimes b \mapsto f_{a\otimes b}$, where $f_{a\otimes b}(1\otimes 1) = af(1 \otimes b)$. The maps $(\cdot u)^*$ and $(\cdot v)^*$ are then the same as the maps $\cdot u$ and $\cdot v$, respectively. Moreover, $B^e$ is free over itself, so $\Hom_B(B^e,B)=(B^e)^*$ is a projective $B^e$-module.  This implies (\ref{eq:Tate10}) is an exact sequence of projective $B^e$-modules; hence, an exact sequence of projective $\D$-modules.

Splicing (\ref{eq:Tate9}) and (\ref{eq:Tate10}) together, we form a $\D$-complete resolution of $B$:
\begin{equation}
\label{eq:Tate11}
\mathbb{P}: \hspace{0.5cm} \cdots \xrightarrow{\cdot u} B^e \xrightarrow{\cdot v} B^e \xrightarrow{\cdot u} B^e \xrightarrow{\xi} B^e \xrightarrow{\cdot u} B^e \xrightarrow{\cdot v} B^e \xrightarrow{\cdot u} \cdots,
\end{equation}
where $\xi = m^* \circ m$. Due to the isomorphism $\HH^*(A) \cong \Ext^*_{\D}(B,A)$, we apply $\Hom_{\D}(-,A)=\widehat{(-)}$ to (\ref{eq:Tate11}):
\begin{equation}
\label{eq:Tate12}
\cdots \xrightarrow{\widehat{\cdot u}} \widehat{B^e} \xrightarrow{\widehat{\cdot v}} \widehat{B^e} \xrightarrow{\widehat{\cdot u}} \widehat{B^e} \xrightarrow{\widehat{\xi}} \widehat{B^e} \xrightarrow{\widehat{\cdot u}} \widehat{B^e} \xrightarrow{\widehat{\cdot v}} \widehat{B^e} \xrightarrow{\widehat{\cdot u}} \cdots, 
\end{equation}
where $\widehat{B^e}$ denotes $\Hom_{\D}(B^e,A)$, and $\widehat{d}(f)=f\circ d$. It is easy to check that the composition of any two consecutive maps $\widehat{d} \circ \widehat{d}$ is equal to $0$, making (\ref{eq:Tate12}) a complex.

For all $f \in \Hom_{B^e}(B^e,A), g \in G$, and $a\otimes b \in B^e$, we have $\Hom_{\D}(B^e,A) \cong \Hom_{B^e}(B^e,A)^G$, where $G$ acts on $\Hom_{B^e}(B^e,A)$ by $(g \cdot f)(a\otimes b) = g \cdot f(g^{-1} \cdot (a\otimes b))$. Note that as a $B^e$-homomorphism, $f$ is completely determined by its value on $1 \otimes 1$. We identify $A$ with $\Hom_\bk(\bk,A) \cong \Hom_{B^e}(B^e,A)$, under the correspondence $a \mapsto f_{a}$, where $f_{a}(1 \otimes 1) = a$, for all $a \in A$. The complex (\ref{eq:Tate12}) becomes: 
\[\cdots \xrightarrow{\widehat{\cdot u}} A^G \xrightarrow{\widehat{\cdot v}} A^G \xrightarrow{\widehat{\cdot u}} A^G \xrightarrow{\widehat{\xi}} A^G \xrightarrow{\widehat{\cdot u}} A^G \xrightarrow{\widehat{\cdot v}} A^G \xrightarrow{\widehat{\cdot u}} \cdots \]
with the actions of $G$ on $A$ depending on the degrees as stated above. The maps $\widehat{\cdot u}$ and $\widehat{\cdot v}$ are:
\begin{align*}
 (\widehat{\cdot u})(\ba) &= (\widehat{\cdot u})f_{\ba}(1 \otimes 1) = f_{\ba}(\cdot u(1 \otimes 1)) = x\ba - \ba x, \\
 (\widehat{\cdot v})(\ba) &= x^{N-1}\ba + x^{N-2}\ba x + x^{N-3}\ba x^2 + \cdots + \ba x^{N-1},
\end{align*}
for all $\ba \in A^G$. We use an analogous argument as in (\cite{BW}, proof of Theorem~2.4) and apply the group actions on $\Hom_{B^e}(B^e,A)$ to take the invariants $A^G$. Since $\chi^{iN}=1$, we find that in 0 and even degrees, $A^G = Z(\kG)$, the center of the group algebra $\kG$, which is $\kG$ itself because $G$ is cyclic. Similarly, in odd degrees, the invariants are spanned by elements of the form $Nxt$, for $t\in G$. However, as $G$ is cyclic generated by $g$, $t=g^j$ for some $j=0,1,\ldots,N-1$, we have $A^G$ is spanned over $\bk$ by $\{x, xg, xg^2, \ldots, xg^{N-1}\}$ in odd degrees. Thus, $\widehat{\cdot v}$ is the $0$-map: $(\widehat{\cdot v})(xg^j) = 0 $, as $x^N = 0$ in $A = B\#\kG$. We then have $\Ker(\widehat{\cdot v})= A^G$ in odd degrees, and $\Imm(\widehat{\cdot v})=0$. Similarly, in even degrees, $\Ker(\widehat{\cdot u}) = \bk$. In odd degrees, $\Imm(\widehat{\cdot u})$ is spanned over $\bk$ by $\{xg, xg^2, \ldots, xg^{N-1}\}$.

We observe that: 
\[\widehat{\xi}: A^G= \text{Span}_{\bk}\{x, xg, xg^2, \ldots, xg^{N-1}\} \rightarrow A^G= \kG\]
maps from degree $-1$ to degree $0$. However, as there is no group element or field element in degree $-1$, $x$ and its powers must be sent to 0. It follows that $\widehat{\xi}$ must be a 0-map. Putting these together, we obtain the Tate-Hochschild cohomology for the Taft algebra $A$, for any integer $n$:
$$\widehat{\HH}^n(A) = \begin{cases}
                      \Ker(\widehat{\cdot u}) / \Imm(\widehat{\xi}) = \bk / 0 = \bk, &\quad n=0 \\
                      \Ker(\widehat{\cdot u}) / \Imm(\widehat{\cdot v}) = \bk/0 = \bk, &\quad n \text{ is even}\\
                      \Ker(\widehat{\cdot v}) / \Imm(\widehat{\cdot u}) = \text{Span}_{\bk}\{x\}, &\quad n \text{ is odd.}
                      \end{cases} $$
Recall the symmetric property that we observed in Example~5.1 for the Sweedler algebra $H_4$:
\[\widehat{\HH}^n(H_4,H_4) \cong \widehat{\HH}^{-(n+1)}(H_4,H_4), \text{ for all } n \in \mathbb{Z}. \]
Indeed, this symmetric result holds in our computation, since any two finite dimensional vector spaces over $\bk$ having the same dimension is isomorphic, $\bk \cong \text{Span}_{\bk}\{x\}$.

\section{Acknowledgement}
The author would like to thank S. Witherspoon for her valuable guidance and support during the preparation of this manuscript. The author is also thankful for the referee's suggestions that reveal results from the literature before which the author was not fully aware of.



\begin{thebibliography}{99}

\bibitem{AvMa} Luchezar L. Avramov and Alex Martsinkovsky, ``Absolute, relative, and Tate cohomology of modules of finite Gorenstein dimension,'' Proc. London Math. Soc., \textbf{85} (2002), no. 3, 393-440. 
\bibitem{BeCa} D. J. Benson and J. F. Carlson, ``Products in negative cohomology,'' J. Pure and Applied Algebra, \textbf{82} (1992), 107-129.
\bibitem{BeJo} P. A. Bergh and D. A. Jorgensen, ``Tate-Hochschild homology and cohomology of Frobenius algebras,'' to appear in J. Noncommutative Geometry, (2013).
\bibitem{Ben1} D. J. Benson, \textit{Representations and Cohomology I: Basic representation theory of finite groups and associative algebras}, Cambridge Studies in Advanced Mathematics, Cambridge Univ. Press, \textbf{30} (1991).
\bibitem{Ben2} D. J. Benson, \textit{Representations and Cohomology II: Cohomology of groups and modules}, Cambridge Studies in Advanced Mathematics, Cambridge Univ. Press, \textbf{31} (1991).
\bibitem{Brown} K. S. Brown, \textit{Cohomology of Groups}, Graduate texts in mathematics, Springer-Verlag, \textbf{87} (1982).
\bibitem{Buch} R.-O. Buchweitz, ``Maximal Cohen-Macaulay modules and Tate-cohomology over Gorenstein rings,'' Preprint, Univ. Hannover, (1986), http://hdl.handle.net/1807/16682.
\bibitem{BW} S. M. Burciu and S. J. Witherspoon, ``Hochschild cohomology of smash products and rank one Hopf algebras,'' Biblioteca de la Revista Matematica Iberoamericana Actas del ``XVI Coloquio Latinoamericano de Algebra,'' (Colonia, Uruguay, 2005), (2007), 153-170. 
\bibitem{CaEi} H. Cartan and S. Eilenberg, \textit{Homological Algebra}, Princeton University Press, Princeton, NJ, (1956).
\bibitem{EiMac} S. Eilenberg and S. Mac Lane, ``Cohomology Theory in Abstract Groups. I,'' Ann. of Math., \textbf{48} (1947), no. 1, 51-78.
\bibitem{EuSc} C.-H. Eu and T. Schedler, `` Calabi-Yau Frobenius algebras,'' J. Algebra, \textbf{321} (2009), no. 3, 774-815.
\bibitem{FiMoSc} D. Fischman, S. Montgomery and H. Schneider, ``Frobenius extensions of subalgebras of Hopf algebras,'' Trans. Amer. Math. Soc., \textbf{349} (1996), 4857-4895.
\bibitem{GiKu} V. Ginzburg and S. Kumar, ``Cohomology of quantum groups at roots of unity,'' Duke Math. J., \textbf{69} (1993), 179-198.
\bibitem{Goi}	F. Goichot, ``Homologie de Tate-Vogel \'{e}quivariante,'' J. Pure and Applied Algebra, \textbf{82} (1992), 39-64.
\bibitem{Mont} S. Montgomery, \textit{Hopf Algebras and Their Actions on Rings}, CBMS Regional Conference Series in Mathematics, Amer. Math. Soc., \textbf{82} (1993).
\bibitem{Mont2} S. Montgomery, ``Representation theory of semisimple Hopf algebras'', \textit{Algebra - Representation Theory}, (Constanta, 2000), K. W. Roggenkamp and M. Stefanescu, editors, NATO Sci. Ser. II Math. Phys. Chem. \textbf{28}, Kluwer, Dordrecht (2001), 189-218.
\bibitem{Nakayama} T. Nakayama, ``On the complete cohomology theory of Frobenius algebras,'' Osaka Math. J., \textbf{9} (1957), 165-187.
\bibitem{PW} J. Pevtsova and S. J. Witherspoon, ``Varieties for modules of quantum elementary abelian groups,'' Algebras and Rep. Th., \textbf{12} (2009), no. 6, 567-595. 
\bibitem{Ste} D. Stefan, ``Hochschild cohomology on Hopf Galois extensions,'' J. Pure and Applied Algebra, \textbf{103} (1995), 221-233.
\bibitem{Tate} J. Tate, ``The higher dimensional cohomology groups of class field theory,''  Ann. of Math., \textbf{56} (1952), no. 2, 294-297.
\bibitem{Weibel} C. Weibel, \textit{An Introduction to Homological Algebra}, Cambridge Univ. Press, (1994).

\end{thebibliography}
\end{document}